\documentclass[11pt]{amsart} 

\usepackage{graphics}
\usepackage{amscd} 

\newtheorem{theorem}{Theorem}

\newtheorem{corollary}[theorem]{Corollary}
\newtheorem{proposition}[theorem]{Proposition}
\newtheorem{ques}[theorem]{Question}
\newtheorem{exam}[theorem]{Example}
\theoremstyle{definition}
\newtheorem{definition}[theorem]{Definition}
\theoremstyle{remark}
\newtheorem{rem}[theorem]{Remark}

\begin{document}

\title[Homeomorphisms of Bagpipes]{Homeomorphisms of Bagpipes}
\author[David Gauld]{David Gauld}
\address{Department of Mathematics, The University of
  Auckland, Private Bag 92019, Auckland, New Zealand}
  \email{d.gauld@auckland.ac.nz}

\thanks{Supported in part by the Marsden Fund Council from Government funding, 
administered by the Royal Society of New Zealand}



\begin{abstract}
We investigate the mapping class group of an orientable $\omega$-bounded surface. Such a surface splits, by Nyikos's Bagpipe Theorem, into a union of a bag (a compact surface with boundary) and finitely many long pipes. The subgroup consisting of classes of homeomorphisms fixing the boundary of the bag is a normal subgroup and is a homomorphic image of the product of mapping class groups of the bag and the pipes.
\end{abstract}
\maketitle

\section{Introduction} \label{Introduction}

A quarter of a century ago Nyikos gave in his Bagpipe Theorem a decomposition of an $
\omega$-bounded surface $F=K\cup\left(\bigcup_{i=1}^nP_i\right)$, where $K$ is a 
compact surface with finitely many holes and $P_1,\dots P_n$ is a collection of long 
pipes attached at these holes. This decomposition makes feasible a partial analysis of 
the mapping class group of such a surface. The mapping class group $\mathcal M(F)$ is 
the quotient of the group of homeomorphisms of $F$ by the normal subgroup of those 
isotopic to the identity. In this paper we will confine our consideration to orientable 
manifolds where convention dictates that the homeomorphisms are orientation-preserving. 

We call a connected Hausdorff space each point of which has a neighbourhood homeomorphic to $\mathbb R^n$ an \emph{$n$-manifold} or just \emph{manifold}. A $2$-manifold will be called a \emph{surface}. We denote the identity homeomorphism by {\bf 1}. The \emph{open long ray} and \emph{long line} will be denoted respectively by $\mathbb L_+$ and $\mathbb L$. A topological space $X$ is called \emph{$\omega$-bounded}, \cite[p. 662]{N},  provided that every countable subset of $X$ has compact closure. The term \emph{long pipe} is defined precisely in Section \ref{Bagpipes and Homeomorphisms} but essentially it is the union of an increasing $\omega_1$-sequence of open sets homeomorphic to $\mathbb S^1\times[0,\infty)$. Recall that two homeomorphisms $h_0,h_1:X\to X$ are \emph{isotopic} provided that there is a continuous function $H:X\times[0,1]\to X$ such that if $H_t:X\to X$ is defined for each $t\in[0,1]$ by $H_t(x)=H(x,t)$ then $H_t$ is a homeomorphism and $H_0=h_0$ and $H_1=h_1$. The map $H$ is called an \emph{isotopy} from $h_0$ to $h_1$. If $h_0$ is isotopic to $h_1$ then we write $h_0\cong h_1$.

Our first main result, Theorem \ref{PreserveBag}, shows that for $F$ as above, any 
homeomorphism $g:F\to F$ is isotopic to a homeomorphism $h$ which leaves the bag 
$K$ invariant and permutes the pipes. Consequently some finite power of $h$, say $h^m$, 
leaves both the bag and each pipe invariant. We show that $h$ could have been chosen 
so that $h^m$ fixes $\d K$ pointwise. Turning to the subgroup $\mathcal N(F)<\mathcal M(F)$ 
consisting of those classes of homeomorphisms leaving bag and pipes invariant, we find in 
Theorem \ref{direct product} that $\mathcal N(F)$ is a normal subgroup and is the 
homomorphic image of the direct product of the groups $\mathcal M(K)$, $\mathcal M(P_1)$, 
\dots, $\mathcal M(P_n)$.

In Section \ref{Type I Spaces and Homeomorphisms} we consider a fairly general 
situation involving homeomorphisms and isotopies of Type I spaces, of which $\omega$-bounded surfaces are a special case. A space $X$ is of \emph{Type I}, \cite[p. 639]{N}, 
provided that it is the union of an $\omega_1$-sequence $\langle U_\alpha\rangle$ of open subsets such that $\overline{U_\alpha}\subset U_\beta$ and $U_\alpha$ is Lindel\"of whenever $\alpha<\beta<\omega_1$. Corollary 5.4 of \cite{N} states that a manifold is $\omega$-bounded if and only if it is countably compact and of Type I.

Section \ref{Bagpipes and Homeomorphisms} applies the results of Section \ref{Type I 
Spaces and Homeomorphisms} to $\omega$-bounded surfaces and it is here that we prove 
Theorem \ref{PreserveBag}. In Section \ref{Mapping Class Groups} we look specifically 
at the mapping class groups and there we verify the facts noted above concerning $\mathcal 
N(F)$.

In Section \ref{Long Pipes} we discuss briefly some possible forms of the mapping class 
group of a long pipe. This is followed by a short section where we raise some questions 
for further investigation.

\section{Type I Spaces and Homeomorphisms} \label{Type I Spaces and 
Homeomorphisms}

\begin{theorem} \label{Type I homeom}
Suppose that $X$ is of Type I, with $X=\cup_{\alpha\in\omega_1}U_\alpha$ as in the definition, and 
that $h:X\to X$ is a homeomorphism. Then $\{\alpha<\omega_1\ /\ h(\cup_{\beta<\alpha}U_\beta)=\cup_
{\beta<\alpha}U_\beta\}$ is a closed unbounded subset of $\omega_1$.
\end{theorem}
\noindent {\bf Proof.} Denote the set by $S$.
\begin{enumerate}
\item {\bf $S$ is closed.} Suppose that $\langle\alpha_m\rangle$ is an increasing sequence of 
elements of $S$ with $\alpha_m\uparrow\alpha$. Then
\begin{eqnarray*} 
h\left(\bigcup{_{\beta<\alpha}}U_\beta\right)& = & h\left(\bigcup{_{m=0}^\infty}_{,\ \beta<\alpha_m}U_
\beta\right)\\ 
& = & \bigcup{_{m=0}^\infty}\ h\left(\bigcup{_{\beta<\alpha_m}}U_\beta\right)\\
 & = & \bigcup{_{m=0}^\infty}\bigcup{_{\beta<\alpha_m}}U_\beta\\
 & = & \bigcup{_{\beta<\alpha}}U_\beta,
\end{eqnarray*}
so $\alpha\in S$.
\item {\bf $S$ is unbounded.} Given any $\alpha_0<\omega_1$, construct an increasing 
sequence $\langle\alpha_m\rangle$ as follows. Using Lindel\"ofness of $U_{\alpha_m}$, given $\alpha_m$ with $m$ even, choose $\alpha_{m+1}>
\alpha_m$ so that
$h(U_{\alpha_m})\subset U_{\alpha_{m+1}}$, 
and given $\alpha_m$ with $m$ odd, choose $\alpha_{m+1}>\alpha_m$ so that
$h^{-1}(U_{\alpha_m})\subset U_{\alpha_{m+1}}$.
As an increasing sequence, $\langle\alpha_m\rangle$ converges, say to $\alpha$. Moreover
\begin{eqnarray*} 
h\left(\bigcup{_{\beta<\alpha}}U_\beta\right) & = & h\left(\bigcup{_{m=0, m\mbox{\footnotesize\ 
even}}^\infty}U_{\alpha_m}\right)\\ 
& = & \bigcup{_{m=0, m\mbox{\footnotesize\ even}}^\infty}\ h(U_{\alpha_m})\\
 & \subset & \bigcup{_{m=0, m\mbox{\footnotesize\ even}}^\infty}U_{\alpha_{m
+1}}\\
 & = & \bigcup{_{\beta<\alpha}}U_\beta
\end{eqnarray*}
and similarly, by letting $m$ range through the odd integers,
\begin{eqnarray*} 
h^{-1}\left(\bigcup{_{\beta<\alpha}}U_\beta\right) & \subset & \bigcup{_{\beta<\alpha}}U_\beta,
\end{eqnarray*}
so $\alpha\in S$. Note also that $\alpha>\alpha_0$. \hspace*{\fill} \rule{1.5mm}{1.5mm}
\end{enumerate} 

\begin{corollary}\label{Isotope Type I}
Suppose that $X$ is of Type I, with $X=\cup_{\alpha\in\omega_1}U_\alpha$ as in the definition, and 
that $h_t:X\to X$, $t\in [0,1]$, is an isotopy of homeomorphisms. Then
\[\left\{\alpha<\omega_1\ /\ h_t\left(\bigcup{_{\beta<\alpha}}U_\beta\right)=\bigcup{_{\beta<\alpha}}U_\beta 
\mbox{ for all } t\right\} \]
is a closed unbounded subset of $\omega_1$.
\end{corollary}
\noindent {\bf Proof.} Choose any countable dense subset $D\subset[0,1]$. By Theorem 
\ref{Type I homeom}, for each $t\in D$ the set
\[ \left\{\alpha<\omega_1\ /\ h_t\left(\bigcup{_{\beta<\alpha}}U_\beta\right)=\bigcup{_{\beta<\alpha}}U_\beta\right
\}\]
is a closed unbounded subset of $\omega_1$. As a countable intersection of closed 
unbounded subsets of $\omega_1$ is closed and unbounded it follows that
\[ \left\{\alpha<\omega_1\ /\ h_t\left(\bigcup{_{\beta<\alpha}}U_\beta\right)=\bigcup{_{\beta<\alpha}}U_\beta 
\mbox{ for all } t\in D\right\}\]
is a closed unbounded subset of $\omega_1$. As $D$ is dense it follows that 
\[ \left\{\alpha<\omega_1\ /\ h_t\left(\bigcup{_{\beta<\alpha}}U_\beta\right)=\bigcup{_{\beta<\alpha}}U_\beta 
\mbox{ for all } t\right\}\]
is a closed unbounded subset of $\omega_1$.\hspace*{\fill} \rule{1.5mm}{1.5mm}

\begin{corollary}\label{Isotope Product}
Suppose that $M$ is a metrisable manifold and that $h_t:M\times\mathbb L_+\to M\times\mathbb L_+$, $t\in [0,1]
$, is an isotopy of homeomorphisms. Then there are $\alpha<\omega_1$ and an isotopy $g_t:M
\to M$ such that $$h_t(\{x\}\times[\alpha,\omega_1))=\{g_t(x)\}\times[\alpha,\omega_1)$$ for each $x\in M$.
\end{corollary}
\noindent {\bf Proof.} Let $D\subset[0,1]$ be a countable dense subset. As a metrisable manifold, $M$ is also separable; say $E\subset M$ is a countable dense subset. For each $t\in[0,1]$ and $x\in M$ consider the function $\theta_{t,x}:\mathbb L_+\to M$ which sends $y\in\mathbb L_+$ to the first coordinate of $h_t(x,y)$. By applying \cite[Lemma 3.4(iii)]{N}, and recalling that a metrisable manifold embeds in euclidean space, see also \cite[Lemma 4.3]{BGG}, one can find $\alpha_{t,x}<\omega_1$ so that $\theta_{t,x}|[\alpha_{t,x},\omega_1)$ is constant. Let $\alpha<\omega_1$ be an upper bound for the countable set $\{\alpha_{t,x}\ /\ t\in D \mbox{ and } x\in E\}$. Thus when $t\in D$ and $x\in E$, $\theta_{t,x}(y)$ is independent of $y\in\mathbb L_+$ when $y\ge\alpha$. As $D$ and $E$ are dense it follows that $\theta_{t,x}(y)$ is independent of $y\in\mathbb L_+$ when $y\ge\alpha$, for all $t\in[0,1]$ and $x\in M$.

Because $M\times\mathbb L_+$ is of Type I, by Corollary \ref{Isotope Type I} we may assume that $
\alpha$ is such that $h_t(M\times(0,\alpha))=M\times(0,\alpha)$ for each $t\in[0,1]$. It follows that $h_t(M\times\{\alpha\})=M\times\{\alpha\}$. Define $g_t:M\to M$ by letting $g_t(x)=\theta_{t,x}(\alpha)$. Then $g_t$ satisfies the requirements. \hspace*{\fill} \rule{1.5mm}{1.5mm}

\bigskip

\section{Bagpipes and Homeomorphisms} \label{Bagpipes and Homeomorphisms}

\begin{definition}\label{long pipe}
An \emph{open long pipe} (long pipe in \cite[page 662]{N}) is a surface $P$ such that 
$P=\cup_{\alpha<\omega_1}U_\alpha$, where each $U_\alpha$ is an open subset of $P$ and is 
homeomorphic to $\mathbb S^1\times\mathbb R$, such that $\overline{U_\alpha}\subset U_\beta$ and that the boundary of each $U_\alpha$ in $U_\beta$ is homeomorphic to $\mathbb S^1$ whenever $\alpha<\beta$. A \emph{(closed) long pipe} is a surface $P$ with boundary $\mathbb S^1$ such that $P=\cup_{\alpha<\omega_1}U_\alpha$, where each $U_\alpha$ is an open subset of $P$ containing the boundary of $P$ and is homeomorphic to $\mathbb S^1\times[0,\infty)$, such that $\overline{U_\alpha}\subset U_\beta$ and that the boundary of each $U_\alpha$ in $U_\beta$ is homeomorphic to $\mathbb S^1$ whenever $\alpha<\beta$.
\end{definition}

We prefer to reserve the term \emph{long pipe} for a closed long pipe as in Definition \ref
{long pipe} rather than Nyikos's terminology because the pipe then contains a boundary 
component which may be used to attach the pipe to the boundary of the bag. Note that 
the equation $P=\cup_{\alpha<\omega_1}U_\alpha$ displays $P$ as a Type I space.

It is easy to show that an $\omega$-bounded long pipe cannot be Lindel\"of.

We will call a decomposition of an $\omega$-bounded surface $F$ into $K\cup\bigg(\bigcup_
{i=1}^nP_i\bigg)$ as in \cite[Theorem 5.14]{N}, but with closed long pipes,  a \emph
{Nyikos decomposition}.

\begin{corollary}\label{IsotopePipes} Suppose the $\omega$-bounded surface $F$ has Nyikos 
decomposition $K\cup\bigg(\bigcup_{i=1}^nP_i\bigg)$ with $P_i=\cup_{\alpha<\omega_1}U_{i,
\alpha}$  a decomposition as in the definition of closed long pipe, and that $h_t:F\to F$, $t\in
[0,1]$, is an isotopy of homeomorphisms. Then\\
${\bigg\{}\alpha<\omega_1\ /\ h_t\left(K\cup\left(\bigcup{_{i=1}^{n}}_{,\ \beta<\alpha}U_{i,\beta}\right)\right)$

\hspace{30mm}$=K\cup\left(\bigcup{_{i=1}^{n}}_{,\ \beta<\alpha}U_{i,\beta}\right) \mbox{ for all } t{\bigg\}}$\\
is a closed unbounded subset of $\omega_1$.
\end{corollary}
\noindent {\bf Proof.} This follows from Corollary \ref{Isotope Type I} because $F$ is of 
Type I. \hspace*{\fill} \rule{1.5mm}{1.5mm}

\begin{rem}
{\rm The definition of long pipe does not assume that $\cup_{\alpha<\lambda}U_\alpha=U_\lambda$ 
when $\lambda>0$ is a limit ordinal (Nyikos calls sequences where this does hold \emph
{canonical sequences}, \cite[Definition 4.3]{N}). So we cannot replace the equation 
defining the closed bounded subset of $\omega_1$ in Corollary \ref{IsotopePipes} by the 
simpler equation $h\left(K\cup\left(\bigcup{_{i=1}^{n}}U_{i,\alpha}\right)\right)=K\cup\left
(\bigcup{_{i=1}^{n}}U_{i,\alpha}\right)$. Example \ref{BadBoundary} shows neither the 
boundary of $\bigcup_{\gamma<\alpha}U_\gamma$  nor the boundary of $P-\overline{\bigcup_{\gamma<
\alpha}U_\gamma}$ need be homeomorphic to $\mathbb S^1$. It is also worth noting that Nyikos points 
out in \cite[p. 669]{N} that if such `bad' boundaries show up over and over again in one 
decomposition then they show up in all decompositions.}
\end{rem}

\begin{exam}\label{BadBoundary}
{\rm Consider the set\\
\hspace*{10mm} $U=\left\{ \left(e^{2\pi it},x\right)\in\mathbb S^1\times\mathbb R\ /\ t=0 \mbox{\rm\ and } x<-1,
\right.$

\hspace{15mm} $\left.\mbox{\rm\ or } 0<t\le\frac{1}{\pi} \mbox{\rm\ and } x<\sin \frac{1}{t}, 
\mbox{\rm\ or } \frac{1}{\pi}\le t<1 \mbox{\rm\ and } x<0 \right\},$\\
an open subset of $\mathbb S^1\times\mathbb R$ which may be expressed as a union $\cup_{n\in\omega}U_n$ 
where $\overline{U_n}\subset U_{n+1}$ and the boundary of each $U_n$ in $U_{n+1}$ 
is homeomorphic to $\mathbb S^1$. A long pipe $P$ may be constructed so that for each limit 
ordinal $\lambda>0$ and each $\alpha<\lambda$ there is a homeomorphism $\varphi:\mathbb S^1\times\mathbb R\to U_
\lambda$ with $\varphi(
\big(\cup_{\alpha<\lambda} U_\alpha\big)=U$. Then the boundaries of $\bigcup_{\gamma<\lambda}U_\gamma$ 
and $P-\overline{\bigcup_{\gamma<\lambda}U_\gamma}$ are both homeomorphic to the boundary of 
$U$ and hence not homeomorphic to $\mathbb S^1$.}
\end{exam}

\begin{corollary}
Suppose that $h:F\to F$ is a homeomorphism of an $\omega$-bounded surface. Then, 
apart from compact subsets, $h$ permutes the long pipes of $F$.
\end{corollary}

Thus if we take an appropriate power of $h$ then, apart from compact subsets, long 
pipes will be preserved set-wise. Hence in order to study the behaviour of 
homeomorphisms of $F$ we need mainly to study homeomorphisms of long pipes. We 
can improve on this if we are happy to work within isotopy classes.

\begin{theorem}\label{PreserveBag}
 Suppose the $\omega$-bounded surface $F$ has Nyikos decomposition $K\cup\bigg
(\bigcup_{i=1}^nP_i\bigg)$ with $P_i=\cup_{\alpha<\omega_1}U_{i,\alpha}$ as in the definition of 
closed long pipe and that $h:F\to F$ is a homeomorphism. Then there is an isotopy of 
homeomorphisms $h_t:F\to F$ so that
 \begin{itemize}
 \item[(1)] $h_0=h$;
 \item[(2)] $h_1(K)=K$.
\end{itemize}
\end{theorem}
\noindent {\bf Proof.} Using Corollary \ref{IsotopePipes} choose any $\alpha<\omega_1$ so that 
$\alpha>0$ and \[h\left(K\cup\left(\bigcup{_{i=1}^{n}}_{,\ \beta<\alpha}U_{i,\beta}\right)\right)=K\cup
\left(\bigcup{_{i=1}^{n}}_{,\ \beta<\alpha}U_{i,\beta}\right).\]
Choose $\gamma,\delta\in\omega_1$ with $\alpha<\gamma<\delta$ so that $h\left(K\cup\left(\bigcup{_{i=1}^
{n}}U_{i,\gamma}\right)\right)\subset K\cup\left(\bigcup{_{i=1}^{n}}U_{i,\delta}\right)$. Note that 
for each $i=1,\dots,n$ there is well-defined $\tilde\imath$ so that $h\left(\overline{U_{i,
\alpha}}-U_{i,\alpha}\right)\subset U_{\tilde\imath,\delta}$, the correspondence $i\mapsto\tilde
\imath$ being a bijection. For each $i=1,\dots,n$ denote the circle $\partial P_i$ by $C_i$ 
and let $D_i\subset K$ be a circle disjoint from but close to $C_i$ so that when $i\not= j$ 
the closed annular region bounded by $C_i\cup D_i$ is disjoint from the closed annular 
region bounded by $C_j\cup D_j$. Again for each $i=1,\dots,n$ let $A_i$ denote the 
closed region bounded by $D_i$ and $\overline{U_{i,\delta}}-U_{i,\delta}$. For each $i$ the 
set $A_{\tilde\imath}$ is an annulus, by the 2-dimensional Annulus Theorem (see \cite[p.
147]{CV} or \cite[p.91]{Mo}, for example) and $h\left(\overline{U_{i,\alpha}}-U_{i,\alpha}\right)$ is 
a circle running once around the interior of the annulus. Hence there is a 
homeomorphism $g_{\tilde\imath}:A_{\tilde\imath}\to A_{\tilde\imath}$ so that $g_{\tilde
\imath}\left(h\left(\overline{U_{i,\alpha}}-U_{i,\alpha}\right)\right)=C_{\tilde\imath}$ and $g_{\tilde
\imath}$ is {\bf 1} on the circles $\left(\overline{U_{\tilde\imath,\delta}}-U_{\tilde\imath.\delta}\right)
\cup D_{\tilde\imath}$. Moreover ${\bf 1}\cong g_i$ by an isotopy $\bar g_t$ which is {\bf 1} on $
\left(\overline{U_{\tilde\imath,\delta}}-U_{\tilde\imath.\delta}\right)\cup D_{\tilde\imath}$. We 
may define $h_t$ to be $h$ except on $\cup_{i=1}^n\mbox{Int}A_i$ and to be $\bar g_th
$ on $A_i$ for each $i$. \hspace*{\fill} \rule{1.5mm}{1.5mm}

\section{Mapping Class Groups}\label{Mapping Class Groups}

Our first result shows that a homeomorphism leaving the bag invariant is isotopic to a 
homeomorphism a power of which is the identity on the boundary of the bag.

\begin{proposition}\label{pipe base identity}
Suppose that $h:K\to K$ is an orientation-preserving homeomorphism of a compact, 
orientable surface with boundary and $m$ is a positive integer so that $h^m$ is invariant 
on each boundary component of $K$. Then there is an isotopy $h_t:K\to K$ so that 
$h_0=h$ and $h_1^m$ is the identity on $\partial K$.
\end{proposition}
\noindent{\bf Proof.} It suffices to consider the case where $h$ cycles the boundary 
components of $K$. Denote the boundary components by $C_1,\dots,C_m$ and assume 
that $h(C_i)=C_{i+1}$ for each $i$ (counted modulo $m$). Let $g_t:C_1\to C_1$ be an 
isotopy such that $g_0=h^m$ and $g_1={\bf 1}$. Firstly define the isotopy $h_t$ on $\partial K$ as follows: $h_t|C_i=h$ when $i<m$ while $h_t|C_m=g_th^{1-m}$. Note that $h_0=h$ while 
$h_1^m$ on $C_i$, for any $i=1,\dots,m$, is \\
$(h_1|C_{i-1})\dots(h_1|C_1)(h_1|C_m)(h_1|C_{m-1})\dots(h_1|C_i)$\\
$=h^{i-1}(g_1h^{1-m})h^{m-i}={\bf 1}$.

Extend $h_t$ over $K$ as follows. By \cite[Theorem 2]{Br}, $C_m$ is collared in $K$, 
more precisely, there is an embedding $e:C_m\times[0,1]\to K$ so that $e(x,1)=x$ for each 
$x\in C_m$. Define the isotopy of embeddings $\varphi_t:K\to K$ by $\varphi_t(e(x,s))=e(x,(1-
\frac{t}{2})s)$ when $(x,s)\in C_m\times[0,1]$ and $\varphi_t$ the identity on $K-e(C_m\times(0,1])
$. Note that $\varphi_t$ is well-defined on $e(C_m\times\{0\})$ and that $\varphi_t(K)=K-e(C_m\times(1-
\frac{t}{2},1])$. Define
\[ h_t(y)=\left\{\begin{array}{lll}
e(h_{2s+t-2)}(x),s) & \mbox{if} & y=e(x,s) \mbox{ and } s\ge1-\frac{t}{2},\\
\varphi_th\varphi_t^{-1}(y) & \mbox{if} & y=e(x,s) \mbox{ and } s\le1-\frac{t}{2}\\
&& \mbox{ or } y\in K-e(C_m\times(0,1]).
\end{array} \right. \]
If $y=e(x,1-\frac{t}{2})$ then\\ 
$\varphi_th\varphi_t^{-1}(y)=\varphi_the(x,1)=\varphi_te(h(x),1)=e(h(x),1-
\frac{t}{2})=e(h(x),0)$\\
$=e(h_0(x),s)=e(h_{2s+t-2}(x),s)$\\ 
so $h_t$ is well-defined. When 
$s=1$, $h_t(e(x,1))=e(h_t(x),1)=h_t(x)$, so $h_t$ really does extend the function $h_t$ 
already defined. It is routine to show that $h_t$ is an isotopy. It is also routine to show 
that $h_0=h$. \hspace*{\fill} \rule{1.5mm}{1.5mm}

\begin{corollary}\label{alignboundary}
Let $F=K\cup\left(\bigcup_{i=1}^nP_i\right)$ be a Nyikos decomposition of an orientable 
$\omega$-bounded surface and $h:F\to F$ an orientation preserving homeomorphism such 
that $h(K)=K$ and $h(P_i)=P_i$ for each $i$. Then $h$ is isotopic to a homeomorphism 
which is the identity on $\partial K$.
\end{corollary}
\noindent{\bf Proof.} If we restrict $h$ to $K$ then we may obtain the required isotopy on 
$K$ from Proposition \ref{pipe base identity} with $m=1$. The same procedure as in the 
proof of Proposition \ref{pipe base identity} enables us to extend the isotopy over each 
pipe $P_i$. \hspace*{\fill} \rule{1.5mm}{1.5mm}

We use the following notation, in which $[h]$ denotes the equivalence class of $h$ under 
isotopy (preserving the boundary if applicable) and $F=K\cup\left(\bigcup_{i=1}^nP_i
\right)$ is a Nyikos decomposition of an orientable $\omega$-bounded surface. Recall that 
we are assuming all homeomorphisms are orientation-preserving.
\begin{itemize}
\item $\mathcal M(F)=\{[h]\ /\ h:F\to F \mbox{ is a homeomorphism}\}$;
\item $\mathcal M(K)=\{[h]\ /\ h:K\to K \mbox{ is a homeomorphism and }$

\hspace{60mm}$ h|\partial K={\bf 1}\}$;
\item $\mathcal M(P_i)=\{[h]\ /\ h:P_i\to P_i \mbox{ is a homeomorphism and }$

\hspace{60mm}$ h|\partial P_i={\bf 1}\}$;
\item $\mathcal N(F)=\{[h]\in\mathcal M(F)\ /\ h(K)=K \mbox{ and }h(P_i)=P_i \mbox{ for each } i\}$.
\end{itemize} 

\begin{theorem}\label{direct product}
Let $F=K\cup\left(\bigcup_{i=1}^nP_i\right)$ be a Nyikos decomposition of an orientable 
$\omega$-bounded surface and let $\mathcal M(F)$ and $\mathcal N(F)$ be as above. Then
\begin{itemize}
\item $\mathcal N(F)$ is a normal subgroup of $\mathcal M(F)$;
\item there is an epimorphism $\theta:\mathcal M(K)\times\mathcal M(P_1)\times\dots\times\mathcal M(P_n)\to\mathcal N(F)
$.
\end{itemize}
\end{theorem}
\noindent{\bf Proof.} Firstly suppose that $g,h:F\to F$ are homeomorphisms so that $[h]
\in\mathcal N(F)$. We must show that $[g^{-1}hg]\in\mathcal N(F)$. By Theorem \ref{PreserveBag} 
we may assume that $g(\partial K)=\partial K$ and by Corollary \ref{alignboundary} we may assume that $h$ is the identity on $\partial K$. Then $g^{-1}hg$ is the identity on $\partial K$ and hence $[g^{-1}hg]\in\mathcal N(F)$.

Secondly, define $\theta$ as follows. Given $([h_0],[h_1],\dots,[h_n])\in\mathcal M(K)\times\mathcal M
(P_1)\times\dots\times\mathcal M(P_n)$ let $h:F\to F$ be the homeomorphism which restricts to $h_0$ 
on $K$ and to $h_i$ on $P_i$ for each $i=1,\dots,n$ and set $\theta([h_0],[h_1],\dots,[h_n])
=[h]$. Then $\theta$ is a homomorphism. That $\theta$ is an epimorphism follows from 
Corollary \ref{alignboundary}. \hspace*{\fill} \rule{1.5mm}{1.5mm}

\begin{rem}
{\rm The kernel of the homomorphism $\theta$ of Theorem \ref{direct product} consists of 
those $(n+1)$-tuples $([h_0],[h_1],\dots,[h_n])$ for which the homeomorphism $h$ 
constructed in the proof above is isotopic to the identity. Set $C_i=\partial P_i$ for each $i$. Then $\partial K=\cup_{i=1}^n C_i$. Each homeomorphism $h_i$ $(i>0$) may be chosen to be the identity except in a collared neighbourhood (in $P_i$) of $C_i$ where $h_i$ acts as some number of Dehn twists. The homeomorphism $h_0$ may then be chosen to be the identity except in a collared neighbourhood (in $K$) of $\cup_{i=1}^n C_i$; in the collared neighbourhood of $C_i$, $h_0$ acts by reversing the Dehn twists which $h_i$ applied on the opposite side of $C_i$.
}
\end{rem}

\begin{rem} 
{\rm The mapping class group $\mathcal M(K)$ is well known. In particular if $K$ has genus $
\gamma$ and has $n$ boundary components then \cite[Theorem 1]{Ge} gives $2\gamma+2n-1$ 
specific generators for this group. Each of these generators is determined by a Dehn 
twist around a closed curve which may circle a boundary component, traverse a handle 
or loop around several handles and/or boundary components. Using these generators 
we may construct interesting homeomorphisms of $\omega$-bounded surfaces.
}
\end{rem}

\begin{exam}
{\rm Replace a small polar cap on the 2-sphere $\mathbb S^2$ by a long pipe. Whereas all 
orientation preserving homeomorphisms of $\mathbb S^2$ are isotopic to the identity it may be 
that not all homeomorphisms of the long pipe are. For example in Proposition \ref{order 
n pipe} below we construct a long pipe supporting a homeomorphism whose mapping 
class (allowing the boundary to move) has order $n$ for any preassigned positive integer 
$n$. As a result we may obtain homeomorphisms, rotations of the remains of the sphere, 
which are not isotopic to the identity. We could add a further long pipe near the opposite 
pole
}
\end{exam}

\begin{exam}
{\rm We may obtain a surface of higher genus and having more long pipes by spreading 
$n$ handles and $n$ mutually homeomorphic long pipes around the equator, arranged 
in such a way that a rotation of the sphere through $\frac{2\pi}{n}$ takes each handle 
and long pipe to the adjacent handle (respectively long pipe). More bands of $n$ handles 
or long pipes or both may be distributed along other lines of latitude. Of course the long 
pipes within a particular band must be mutually homeomorphic. See figure \ref{Bristly 
Sphere}.
}
\end{exam}

\begin{figure}
\scalebox{0.5}{\includegraphics{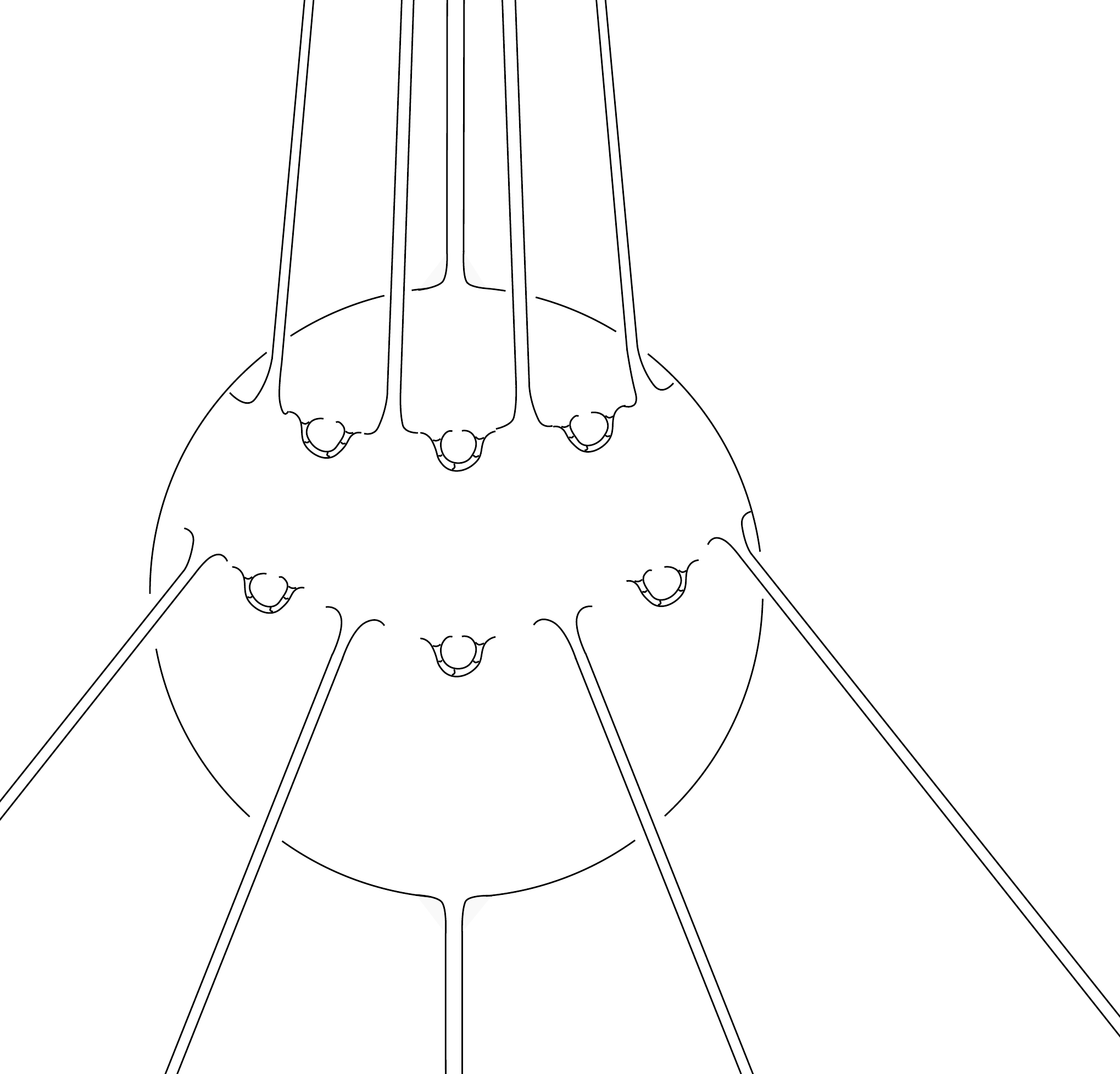}}
\caption{A bristly sphere.}\label{Bristly Sphere}
\end{figure}

\begin{exam}
{\rm An interesting question is whether there are homeomorphisms a finite power of 
which are isotopic to the identity: thus the mapping class group has torsion. In the case 
of the torus $\mathbb T^2$ we may look at the group $GL(2,\mathbb Z)$ as each isotopy class of 
homeomorphisms of $\mathbb T^2$ may be represented by an element of this group. If an 
element $A$ of $GL(2,\mathbb Z)$ is of finite order then its eigenvalues must be roots of the 
cyclotomic polynomial of degree $n$ which must therefore divide the characteristic 
polynomial of $A$ and hence have degree 2. It follows that $n$ can be only 1, 2, 3, 4 or 
6. In \cite[Theorem 2]{Me} it is shown that if $A\in GL(2,\mathbb Z)$ is periodic then $A$ is 
conjugate to one of the six matrices
\begin{itemize}
\item $\left[\begin{array}{rr} 0 & 1\\ 1 & 0 \end{array}\right]$, which has order 2, switching 
the coordinates;
\item $\left[\begin{array}{rr} -1 & 0\\ 0 & 1 \end{array}\right] 
=\left[\begin{array}{rr} 0 & 1\\ 1 & 0 \end{array}\right]\left[\begin{array}{rr} 0 & 1\\ -1 & 0 
\end{array}\right]$, which has order 2, reversing one coordinate;
\item $\left[\begin{array}{rr} -1 & 0\\ 0 & -1 \end{array}\right] 
=\left[\begin{array}{rr} 0 & 1\\ -1 & 0 \end{array}\right]^2$, which has order 2, rotating 
through 180$^\circ$;
\item $\left[\begin{array}{rr} -1 & 1\\ -1 & 0 \end{array}\right] 
=\left[\begin{array}{rr} 0 & 1\\ -1 & 1 \end{array}\right]^2$, which has order 3;
\item $\left[\begin{array}{rr} 0 & 1\\ -1 & 0 \end{array}\right]$, which has order 4, rotating 
through 90$^\circ$;
\item $\left[\begin{array}{rr} 0 & 1\\ -1 & 1 \end{array}\right]$, which has order 6.
\end{itemize}

Note that some of these matrices have determinant $-1$ and hence give rise to 
orientation reversing homeomorphisms. As in the previous example we may add 
strategically placed long pipes and handles to add to the interest.
}
\end{exam}

\section{Long Pipes}\label{Long Pipes}

It is almost hopeless to attempt to determine the mapping class group of all long pipes in 
the same way as the mapping class group has been determined for compact surfaces 
with boundary as in \cite{Ge}. Indeed, as noted in \cite[p. 669]{N}, there are $2^
{\aleph_1}$ topologically distinct long pipes. In this section we describe a few of these, 
especially one possessing homeomorphisms of order $n$ for a fixed given $n\in\mathbb N$.

The simplest example of a long pipe is $\mathbb S^1\times\mathbb L_+$, the product of a circle with the 
closed long ray. Any orientation preserving homeomorphism is isotopic to the identity. 

Before introducing the next example we introduce another concept. Say that a 
homeomorphism $h:X\to X$ has \emph{isotopy order $n$} provided that $n$ is the least 
positive integer for which $h^n={\bf 1}$.

The long plane $\mathbb L^2$ with an open disc removed is another reasonably simple example 
of a long pipe. The mapping class group of $\mathbb L^2$ has been determined in \cite[Theorem 
1.3]{BDG} where it is seen that the isotopy order of any homeomorphism of $\mathbb L^2$ (and 
hence of the long pipe obtained from it) is 1, 2 or 4. This long pipe inspires the next 
proposition.

\begin{proposition} \label{order n pipe}
For every natural number $n$ there is a long pipe $P$ and a homeomorphism $h:P\to P
$ whose (isotopy) order is $n$.
\end{proposition}
\noindent {\bf Proof.} Given $n$ construct the long pipe $P$ as follows. Take $n$ copies 
of the truncated first octant $\{(x,y)\in\mathbb L^2\ /\ x\ge y\ge0 \mbox{ and } x\ge 1\}$, say $
\mathbb O_1,\dots,\mathbb O_n$; denote the point of $\mathbb O_i$ corresponding to $(x,y)$ in the first octant 
by $(x,y)_i$. Denote by $\mathbb O_i^0$ (resp. $\mathbb O_i^1$) the edge of $\mathbb O_i$ corresponding to 
the line $y=0$ (resp. $y=x$): as noted in \cite[Example 3.8]{N} these two subsets behave 
very differently in $\mathbb O_i$ even though each is homeomorphic to the closed long ray. Set 
$P=\big(\dot{\cup}_{i=1}^n\mathbb O_i\big)/\sim$, where $\sim$ identifies $(x,x)_i\in\mathbb O_i^1$ with 
$(x,0)_{i+1}\in\mathbb O_{i+1}^0$ for any $i=1,\dots, n$, where the subscript $i+1$ is modulo $n
$. Define $h:P\to P$ by $h((x,y)_i)=(x,y)_{i+1}$. \hspace*{\fill} \rule{1.5mm}{1.5mm}

Proposition \ref{order n pipe} tells us that if there is a homeomorphism of a compact 
surface with one boundary component which has (isotopy) order $n$ then there is a 
homeomorphism of a surface with one long pipe having (isotopy) order $n$. As an 
example take a surface of genus $n$ where the $n$ handles are spread symmetrically 
about a central axis like the petals of a simple flower and a small disc centred at one 
point where the axis cuts the surface is removed.

On the other hand there are long pipes where no homeomorphism has finite isotopy 
order. The simplest example is $\mathbb S^1\times\mathbb L_+$ but $\mathcal M(\mathbb S^1\times\mathbb L_+)$ has a single 
element so perhaps is uninteresting. A more interesting example is obtained as in 
Proposition \ref{order n pipe} with $n=1$; in this case any homeomorphism $h:P\to P$ is 
fixed (up to isotopy) outside a bounded set so $[h]$ is determined by a number of Dehn 
twists on an annulus bounded on one side by $\partial P$ and hence $\mathcal M(P)\approx\mathbb Z$.

\section{Questions}\label{Questions}

\begin{ques}
Is there a long pipe where no homeomorphism has finite isotopy order?
\end{ques}
Of course in this question we need to be careful what isotopies we allow. The intention 
here is that isotopies need not be the identity on the boundary.

\begin{ques}
Suppose that $[h]\in\mathcal M(P)$ has finite order, where $P$ is a long pipe. Does it follow 
that $h$ is isotopic to the identity? In other words is it true that $\mathcal M(P)$ has no 
torsion?
\end{ques}

\end{document}